\documentclass[12pt]{article}
\usepackage{amsmath, amsfonts, amssymb, amsthm} 
\newtheorem{theorem}{Theorem}[section]
\newtheorem{proposition}[theorem]{Proposition}
\newtheorem{definition}[theorem]{Definition}
\newtheorem{remark}[theorem]{Remark}

\newtheorem{corollary}[theorem]{Corollary}
\newtheorem{lemma}[theorem]{Lemma}

\def\R{\mathbb R} \def\Z{\mathbb Z} \def\C{\mathbb C} 

\def\N{\mathbb N}
\def\C{{\mathbb C}} 
\def\Q{\mathbb Q}

\def\T{\mathbb T}

\def\Ad{\hbox{\rm Ad}\,}

\def\ad{\mathrm {ad}}

\def\trace{{\rm trace}}

\def\<{\,<\!}
\def\>{\!>\,}
\usepackage{graphics}
\usepackage{epsf}
\usepackage{epsfig}
\begin{document}

\title{On the radicals of exponential Lie groups} 

\author{S.G. Dani}
\maketitle
\begin{abstract}
Let $G$ be a connected exponential Lie group and $R$ be the solvable radical
of $G$. We describe a condition on $G/R$ under which one can then conclude
that $R$ is an exponential Lie group. The condition holds in particular when
$G$ is a complex Lie group and this yilds a stronger version of a result of
Moskowitz and Sacksteder \cite{MS} on the center of a complex exponential Lie
group being connected. Along the way we prove a criterion for elements
from certain subsets of a solvable Lie group to be exponential, which would
be of independent interest.
\end{abstract}


\section{Introduction}

Let $G$ be a connected Lie group,  $\frak G$ be the Lie algebra of $G$,
and let $\exp: \frak G \to G$ be the exponential map. 
An element $g\in G$ is said to be
{\it exponential} if there exists $X\in \frak G$ such that $g=\exp X$, namely
if $g$ is contained in the image of the exponential map. The Lie group $G$ is said to be exponential if every element is exponential, or equivaently if 
the exponential map is surjective. There has been considerable interest
in literature to understand which Lie groups are exponential. On the whole
satisfactory results are known when $G$ is either solvable or semisimple
(see \cite{DH} for the status until the mid-1990's and  \cite{W}, \cite{W2} and \cite{W3} 
 for later work). A general Lie group is an almost 
semidirect product of a semisimple and a solvable subgroup (namely 
a semisimple Levi subgroup and the solvable radical respectively),
but not much clarity is attained so far as to which semidirect products 
are exponential Lie groups. Some results in this respect may be found in 
\cite{DM}, \cite{Dj}, \cite{MS} and  \cite{C}.

In this context we consider the following question. 
 Given an exponential Lie group $G$ when can we 
conclude that its radical $R$ is exponential? 
We describe a condition on  $G/R$,  the semisimple quotient Lie group, 
which enables such a conclusion (see Theorem~\ref{thm});
 the condition involved holds for a large class
 of semisimple Lie groups (see Theorem~\ref{ss}). 
 
The question as above arose in the context of the following. 
The main theorem (Theorem~1) in  \cite{MS} states  that a 
connected complex Lie group $G$ is exponential if and only if its center is 
connected and the adjoint group $\Ad (G)$ is exponential; as has been noted
in the review of the paper in Mathematical Reviews the proof of the theorem
given in the paper is valid only under an additional condition that $G$ 
is ``semi-algebraic'' (we shall not go into the details of the condition
as it does not concern us in the sequel). The ``if'' part, 
and also part of the converse, is relatively easy to see, so the main content 
of such a statement (when valid) is that the center of a connected complex
exponential Lie group is connected; moreover,  
a complex semisimple exponential Lie group has trivial center, and hence conclusion is
equivalent to that the center of the radical is connected. 
This would follow directly, if we conclude that the radical is exponential, 
since the center of an exponential solvable Lie group is always connected (see 
in particular, Corollary~\ref{cor}, infra). 
Our results apply to this case and in particular yield
the above mentioned statement from \cite{MS} as a corollary in a special 
case. Our technique 
however is independent of~\cite{MS}. 

We now formulate the condition involved and state the main results. 

To begin, we recall that a one-parameter subgroup of a Lie group $L$ is 
said to be 
{\it unipotent} if it is of the form $\{\exp tX\}_{t\in \R}$, where $X$ 
is an element of the Lie algebra of $L$ such that the adjoint transformation 
$\ad X$ is a nilpotent linear tranformation of the Lie algebra. 

\begin{definition}{\rm We say that a connected Lie group $L$ satisfies 
{\it condition $\frak U$} if there exists a unipotent one-parameter subgroup
$U$ of $L$ such that the centralizer of $U$ in $L$, viz. $\{x\in L\mid xu=ux
\hbox{ \rm for all } u\in U\}$, does not contain any compact subgroup of 
positive dimension. } 

\end{definition}



\begin{theorem}\label{thm} 
Let $G$ be an exponential  Lie group and $R$ be the radical of $G$. 
Let $S=G/R$ and suppose that $S$ satisfies 
condition $\frak U$. 
Then $R$ is an exponential Lie group. 

\end{theorem}

In the context of the theorem it would be of interest to know when a 
semisimple Lie group satisfies condition $\frak U$.  
Let $S$ be a semisimple Lie group and let $S=KAN$ be an Iwasawa 
decomposition of $S$, namely $K, A$ and $N$ are closed connected 
subgroups of $S$, and if $S^*$ is the adjoint group of $S$ and $\Ad :S\to S^*$
is the adjoint representatiton, then $\Ad (K)$ is a maximal compact subgroup
of $S^*$, $\Ad (A)$ is a maximal
connected subgroup whose action on the Lie algebra is diagonalisable  
over $\R$, and 
$N$ is a simply connected nilpotent Lie subgroup of $S$ normalised by 
$A$ (see for example \cite{Hel}); we note that $K$ contains the center 
of $S$ and $K$  itself is not compact when the center is infinite, as in the case of 
some nonlinear semisimple Lie groups. Let $M$ be the centraliser of $A$ in 
$K$. Then $M$ is 
a closed (not necessarily connected) subgroup of $S$. Also $M$  normalises 
 $N$ and  $P=MAN$ is 
a minimal parabolic subgroup of $S$. 

With regard to condition $\frak U$ we prove the
following. 

\begin{theorem}\label{ss}
Let the notation be as above and suppose that $M$ is abelian. Then $S$ satisfies condition 
$\frak U$. In particular if $S$ 
is a $\R$-split semisimple Lie group or a complex semisimple Lie group then 
it satisfies condition $\frak U$. 

\end{theorem}

\begin{remark}
{\rm It has also been possible to prove (without reference to $M$ being 
abelian) that if $N$ contains an 
$\R$-regular unipotent element of $S$, in 
the sense of \cite{A}, which is not centralised by 
 any compact subgroup of $P$ of positive dimension then $S$ satisfies 
condition $\frak U$. The condition holds when $M$ is abelian. 
 It has however not been possible to ascertain whether there are semisimple Lie groups 
 $S$ with $M$ nonabelian for which this holds. We shall therefore not go 
into the technical details
 of this possible generalisation, at this stage.

}
\end{remark}
When $S$ is the group of $\R$-elements of a semisimple algebraic group 
$\underline S$ defined
over $\R$ (or the connected component of the identity in such a group) the 
condition of $M$ being abelian is equivalent to $\underline S$ being a quasi-split
group, namely one admitting a Borel subgroup defined over $\R$ (see \cite{R} 
and \cite{C}). Thus 
Theorem~\ref{ss}
applies in this case; we note however that $S$ as in Theorem~\ref{ss} need not 
be linear. 


Theorems~\ref{thm} and \ref{ss} together imply the following, which yields 
in particular the statement  from \cite{MS} recalled above.  

\begin{corollary}
Let $G$ be an exponential Lie group,  $R$ be the radical of $G$ and 
$S=G/R$. Suppose  $S$ satisfies the condition as in Theorem~\ref{ss}.
 Then $R$ is an exponential Lie group. In particular
the center of $R$ is connected. 

\end{corollary}

We note that for condition $\frak U$ to hold for a semisimple Lie group 
it has to be noncompact. Also, not all noncompact semisimple Lie 
groups satisfy condition $\frak U$; this can be readily verified for the group 
$SU(n,1)$ 
for instance, which is in fact one of the exponential Lie groups 
(see \cite{DN},\cite{W2}).  The question 
whether a Lie group $G$ being exponential implies that its radical is
exponential remains open in the cases when $G/R$ is one of these 
semisimple Lie groups. 

Towards proving Theorem~\ref{thm} we also prove a result on 
 exponential elements in solvable 
Lie groups, which may be of independent interest.  
The question of exponentiality of solvable Lie groups was considered earlier
in \cite{W} (and earlier papers cited there), but the focus there has been on 
criteria for the group to be 
exponential, namely {\it all} elements being exponential. On the 
other hand  our result 
below, Theorem~\ref{solvable}, is about exponentiality of elements from 
certain subsets. Moreover, 
the argument in \cite{W} involves the theory of Cartan subgroups, which is 
technically more intricate; our argument is based on more elementary 
considerations. 

The paper is organised as follows. The result on exponential elements
 in solvable Lie groups,  
Theorem~\ref{solvable},  will be taken 
up in \S 2, and in \S 3 the results are applied to prove Theorem~\ref{thm}
and some related results. In \S 4 we discuss when condition $\frak U$ is 
satisfied, and prove Theorem~\ref{ss}.

\section{Exponential elements in solvable Lie groups}

Let $H$ be a connected solvable Lie group and $N$ be a nilpotent simply 
connected
closed normal Lie subgroup of $H$ such that $H/N$ is abelian.  We denote 
by $\frak N (H)$ 
the class of closed connected normal subgroups of $H$ contained in $N$. 
We shall be considering pairs of the form 
$(M,M')$, with $M,M' \in \frak N (H)$, and $M' \subset M$, and 
 for brevity we shall refer to such a pair as a {\it subquotient} (of $N$, with respect
 to the action of $H$). 
A subquetient $(M,M')$ is said to be abelian if $M/M'$ is abelian; 
as 
$N$ is simply connected, in this case $M/M'$ is a vector space over $\R$. 
The conjugation action of $H$ on itself induces an action on each $M\in \frak N(H)$, by restriction,
and hence on  all subquotients of $H$.  
An abelian subquotient $(M,M')$,  $M, M'\in \frak N(H)$, will be called an {\it irreducible
subquotient}  if  the $H$-action on $M/M'$ is irreducible.  We 
note that for any irreducible subquotient $(M,M')$, $M,M'\in \frak N(H)$, the action of $N$
on $M/M'$ is trivial, and hence the action factors to an irreducible linear action of 
$A$ on $M/M'$. 


\begin{definition}{\rm Let the notation be as above. Let $a\in A=H/N$ and 
$B$ be a one-parameter subgroup of $A$ containing $a$. 
We say that the pair $(a,B)$ is of {\it type~$\cal E$} if the following
holds: if for an irreducible 
subquotient  $(M,M')$, $M,M'\in \frak N(H)$,  
the action of $a$ on $M/M'$ is trivial, then the action of $B$ on $M/M'$ 
is trivial; the element $a$ is said to be type $\cal E$ if there exists 
a one-parameter 
subgroup  $B$  of $A$  such that the pair $(a,B)$ is of type $\cal E$. 
}
\end{definition}

We prove the following.

\begin{theorem}\label{solvable}
Let $H$ be a connected solvable Lie group and $N$ be a 
nilpotent simply connected
closed normal Lie subgroup of $H$ such that $H/N$ is abelian. 
Let $A=H/N$. Let $x\in H$ and $a=xN\in A$. Then
the following conditions are equivalent: 

i) $a$ is of type $\cal E$. 

ii) for every $n\in N$, the element  $xn$ is exponential in $H$. 

\end{theorem}

We first consider the following special case, in which we prove some more precise
 results. By a {\it vector subgroup} we mean a subgroup which is topologically 
 isomorphic to $\R^d$ for some $d\geq 1$ (with respect to the induced 
topology). 
A vector subgroup
will be considered equipped with its canonical structure as a vector space over $\R$. 

\begin{proposition}\label{prep}
Let $L$ be a connected Lie group.  Let $V$ be a vector subgroup of $L$
and $P$ be a one-parameter subgroup normalizing $V$,
and such that the conjugation action of $P$ on $V$ is 
irreducible and nontrivial. Let $H=PV$. Let 
$p \in P$ be nontrivial and let $\sigma (p):V\to V$ denote the 
conjugation action of $p$ on $V$. 
Then the following holds: 


i) Every one-parameter subgroup of $L$ contained in $H$ is either
contained in $V$ or is of the form 
 $ wP  w^{-1}$ for some $w\in V$.

ii)  If  $\sigma (p)$ is nontrivial then for any $v\in V$, $pv$ 
is contained in a unique  one-parameter 
subgroup contained in  $H$.
 
iii)  If   $\sigma (p)$ is trivial, then for $v\in V$, 
$v\neq 0$ (the zero element in $V$), $p v$ is not contained 
in any one-parameter subgroup of $H$.

\end{proposition} 

\proof At the outset we note that the subgroup $H$ as in the hypothesis
need be a closed subgroup of $L$. It is however a Lie subgroup whose 
Lie algebra is the sum of the Lie subalgebras of $P$ (which is
one-dimensional) and $V$. Let $\frak H$ be the Lie algebra of $H$; 
we shall realize it as $\R \xi \oplus V$, where $\xi$ is a generator
of the Lie algebra of $P$ and the vector space $V$ is identified with 
its Lie algebra. We now prove the statements as in the Proposition.

i) For $w\in V$ the Lie subalgebra of $\frak H$ 
corresponding to the one-parameter subgroup $wP w^{-1}$ is spanned by 
$\xi +(\theta (w)-w)$, where $\theta :V\to V$ is given by $\theta (u)=\frac 
{dp_{-t}(u)}{dt}$ for all $u\in V$. We see that as the $P$-action on $V $ 
is nontrivial and irreducible the map $w\mapsto \theta (w)-w$ 
is surjective, and hence all elements of the form $\xi + v$, $v\in V$, are 
among the tangents to  $wP w^{-1}$, $w\in V$. Since scalar multiples of 
these cover all elements that are not contained in $V$, it follows that 
every one-parameter subgroup which is not contained in $V$ is of the form 
$wP w^{-1}$  (upto scaling of the parameter). This proves  (i). 
 
 ii) Let $p\in P$ be such that  $\sigma (p)$ is nontrivial. As the $P$-action 
 on $V$ is irreducible and nontrivial  this
implies in particular that  $1$ is not an eigenvalue of $\sigma (p)$. Hence
$\sigma (p)^{-1} -I$, where $I$ is the identity transformation, 
is invertible. Now let  $v\in V$ be given. Then 
there exists $w\in V$ such that $v=\sigma (p) ^{-1}(w)-w$. 
Hence, in $H$, $v =(p^{-1} wp)w^{-1}$. Thus   $p v 
=wp w^{-1}$, and it is contained in the one-parameter subgroup $wP
w^{-1}$.  

We  note that if $pv$ is contained 
in $wP w^{-1}$ for some $w\in V$ then in fact $p v=wp w^{-1}$; if 
$q\in P$ is such that $pv=wqw^{-1}=q(q^{-1}wqw^{-1})\in qV$, then $p^{-1}q\in V$
and hence $pv=wqw^{-1}=(wpw^{-1})(wp^{-1}qw^{-1})= wpw^{-1}$. 
Now if $p v$  is contained in 
$w_1P w_1^{-1}$ and $w_2P w_2^{-1}$, $w_1,w_2 \in V$,
we have $w_1p w_1^{-1}=w_2p w_2^{-1}$, which means that $w_2^{-1}w_1$
(or $w_1-w_2$ is addivive notation) is fixed by $\sigma (p)$. 
Since $\sigma (p)$ is nontrivial and the $P$-action is irreducible 
this implies that $w_1=w_2$. Thus 
one-parameter subgroup containing $pv$ is 
unique. This proves~(ii).

iii) Now let $p\in P$ be such that  $\sigma (p)$ is trivial.
Let $v\in V$, $v\neq 0$ and 
suppose that $pv$ is contained in $wP w^{-1}$ for 
some $w\in V$. Then, as before, we have $p v=wp w^{-1}$. 
Hence $v=p^{-1}wp w^{-1}=
\sigma (p)^{-1}(w)-w=0$, contradicting the assumption that $v\neq 0$. Hence   $p v$ 
is not contained in any 
one-parameter subgroup of $H$. This proves~(iii). \qed

In the sequel it will be convenient to use the following terminology. 

\begin{definition}
{\rm Given a Lie group $L$, a 
closed normal subgroup $M$ of $L$ with  $\eta :L\to L/M$ the canonical
quotient map, and 
a one-parameter subgroups  $\Psi$ of $L/M$, by a {\it lift} of $\Psi$
in $L$ we mean  
a one-parameter subgroup $\Phi$ of $L$ such that  $\eta (\Phi) = \Psi$. 
}
\end{definition}

\begin{proposition}\label{countable}
Let the notations $H,N, A, x$ and $a$ be as in Theorem~\ref{solvable}. 
Let $B$ be a one-parameter subgroup of $A$ containing  
$a $, such that $(a,B)$ is of type~$\cal E$. Then for any $n\in N$ there 
exists a lift of $B$ containing $xn$. Moreover, the collection of 
such lifts is countable. 
\end{proposition}

\proof 
We shall show that if $(M,M')$, where $M,M'\in \frak N(H)$, is an 
irreducible subquotient then the following holds:
given  a one-parameter subgroup $\Psi$ of $H/M$ which is a lift of $B$
 and $y\in H$ such that $yN=a$,
 there exist at least one and at most countably many lifts  $\Phi$
 of $\Psi$ in $H/M'$  containing $yM'$. 
Applying this to successive pairs from a sequence  $N_0, N_1, \dots , N_k$ 
 in $\frak N(H)$ such that $N_0=N$,
$N_k$ is trivial, and for all $j=0,\dots , k-1$, $N_{j+1} \subset N_j$ and 
$(N_j,N_{j+1})$ is an irreducible subquotient, yields 
the statement as in the proposition.  

Now let $(M,M')$ as above, a one-parameter subgroup $\Psi$ of $H/M$ which is a lift of $B$, and 
$y\in H$ such that $yM=a$ be given.  
Let  $P$ be a lift of $\Psi$ in $H/M'$.  We now apply 
Proposition~\ref{prep} to $P(M/M')$ (the product of $P$ and $M/M'$ in $H/M'$), 
with $v=yM'$.  We note that any one-parameter
subgroup of  $H/M'$ which is a lift of $\Psi$ is contained in  $P(M/M')$, 
thus it suffices
to show that  there exists a lift $\Phi$ of $\Psi$ in $P(M/M')$ 
 containing $yM'$.  If the action of
$a$ on $M/M'$ is nontrivial then this is assured by Proposition~\ref{prep}(ii).
 On the other hand if the action of $a$ on $M/M'$ is  trivial
then by hypothesis the action of $B$ is also trivial, and since $P$ is a lift
of $B$ this implies that $P(M/M')$ is an abelian Lie group,  in 
this case the assertion is obvious.  This proves the Proposition. \qed

\begin{proposition}\label{one-step}
Let $H,N, A, x$ and $a=xN$ be as in Theorem~\ref{solvable} and let $B$ be a 
one-parameter subgroup of $A$ containing $a$. Let $M\in \frak N(H)$ be 
an abelian subgroup such that the $A$-action on $M$ is irreducible and 
the following holds:

i)  for the group $H/M$, with $A$ being viewed canonically as $(H/M)/(N/M)$,
the pair $(a,B)$ is of type $\cal E$. 
 
ii) on $M$ the action of $a$ is trivial but the action of $B$ is not trivial. 

\noindent Let $Q$ be a lift of $B$ in $H/M$ and $y\in H$ be such that $yN=a$.  
Then there exist a unique $m\in M$ 
such that $ym$ is contained in a lift of $Q$. 

\end{proposition}

\proof Let  
$P$ be a lift of $Q$ in $H$. Let $p\in P$ be such that $pM=yM$, and 
$m\in M$ such that $p=ym$. Then $pN=yN=a$. 
By condition~(ii) in the hypothesis and Proposition~\ref{prep}(iii), 
applied to the subgroup $PM$,
we get that $pm'$ is not contained in a lift of $Q$ for any nontrivial 
element $m'$ of $M$. 
Thus $m$ is the only element of $M$ such that $ym$ is contained 
in a lift of $Q$. \qed

\begin{proposition}\label{notE}
Let $H,N, A, x$ and $a=xN$ be as in Theorem~\ref{solvable} and let $B$ be a 
one-parameter subgroup of $A$ containing $a$. Suppose that $(a,B)$ is not of
type $\cal E$. Let $E=\{n\in N\mid xn \hbox{ \rm is contained in a lift of } 
B \hbox{ \rm in } H\}$. Then 
$E$ has $0$ Haar measure in $N$.  

\end{proposition}

\proof At the outset we note that the set of exponential elements in a Lie group 
is a Borel subset, and hence it follows that $E$ as in the hypothesis is a Borel subset of $N$. 
As  $(a,B)$ is not of
type $\cal E$ there exists an irreducible subquotient $(M,M')$, 
$M,M'\in \frak N (H)$, such that the $a$-action on $M/M'$ is trivial and the 
$B$ action is nontrivial, and we may without loss of generality assume $M$ to be 
of maximal possible dimension among such pairs. To prove the propositin it
suffices to show that $EM'/M'$ has zero Haar measure in $N/M'$, and hence 
passing to $H/M'$,  we may without loss of generality assume $M'$ to be trivial. Thus $M$ 
is a normal vector subgroup of $H$.  By the maximality 
condition on $M$, $(a,B)$ is of type $\cal E$ for $H/M$. Hence by 
Proposition~\ref{countable} for any $n\in N$ there exist only countably 
many lifts of $B$ in $H/M$ containing $xnM$. Let $n\in N$ and $Q$ be 
any lift of $B$ in $H/M$ containing $xnM$. Let $P$ be a lift of $Q$ in $H$
and $m\in M$ be such that $xnm\in P$. 
We now apply 
Proposition~\ref{prep}(iii), with $M$ in place of $V$ there.  The choice 
of $M$ as above shows that the conjugation action of $xnm$ on $M$ is trivial but the 
action of $P$ is not trivial. By Proposition~\ref{prep}(iii) therefore $xnm$
is the only element in $xnM$ which is contained in a lift of $Q$. Since 
there are only countably many lifts $Q$ of $B$ containing $xnM$, this 
shows that there are only contably many $m$ in $M$ such that $xnm$ 
is contained in a lift of $B$ in $H$. In other words, $E$ as in the 
hypothesis intersects each coset $xnM$, where $n\in N$, in only countably many points. It follows 
that its Haar measure of $E$ must be $0$. \hfill 

\medskip
\noindent{\it Proof of Theorem~\ref{solvable}}:   
Statement (i) $\implies$ (ii) follows immediately from (the existence 
statement in) 
Proposition~\ref{countable}. We now prove the converse. We suppose 
that condition~(ii) holds, but $a$ is not of type~$\cal E$, and arrive
at a contradiction. 

Let $\cal C$ be the class of $L\in \frak N (H)$, such that $xN/L$ 
is not of type $\cal E$ for the group $H/L$. Then $\cal C$ is nonempty 
since by assumption $a$ is not of type~$\cal E$ for $H$, so 
$\cal C$  contains the trivial subgroup. Let $L$ be an 
element of $\cal C$ with maximum possible dimension; we note that $L$ 
is a proper subgroup of $N$, since  the action of $A$ on $N/L$
has to be nontrivial. 
For notational convenience we shall consider $H/L$ as $H$, and respectively
$N/L$ as $N$, and $xL$ as $x$. Then in the modified notation we have that 
$xn$ is exponential in $H$ for all $n\in N$, $xN$ is not of type $\cal E$, 
but $xN/M$ is of type $\cal E$ for every $M\in \frak N(H)$ of positive 
dimension.
We fix a $M\in \frak N(H)$ of positive dimension such that the action 
of $A=H/N$ on $M$ is irreducible.  Let $a=xN$. Then the action of $a$ on 
$M$ is trivial since otherwise $a$ would be of type~$\cal E$ under the 
conditions as above.  

By hypothesis, for each $n\in 
N$ there exists a one-parameter subgroup containing $xn$. 
We note that the set 
of one-parameter subgroups of $A$ containing $a=xN$ is countable, say
$\{B_j\}$ with $j$ running over a countable set; it can be a singleton
set as would be the case when $A$ is simply connected. 

By Proposition~\ref{notE}, applied to $H/M$ there exists $n\in N$ such 
that $xnM$ is not
contained in a lift of $B_j$ for any $j$ such that $(a,B_j)$ is not of 
type $\cal E$ for $H/M$; indeed the set of such $n$ is a set of full 
Haar measure. We fix such $n$ and consider the elements $xnm$, 
$m\in M$. By hypothesis each of them is exponential and hence is contained
in a lift of $B_j$ for some $j$, and by the  choice of $n$, the 
$j$ must be such that $(a,B_j)$ is of type $\cal E$ for $N/M$. Hence to prove
the theorem it suffices
to show that the action of one of these on $M$ is trivial. Consider
any $B_j$ such that $(a,B_j)$ is of type $\cal E$ for $H/M$. Then  
by Proposition~\ref{countable} 
 for any $j$  there exist only countably many lifts
of $B_j$ to $H/M$ containing $xnM$. Let $Q$ be any lift 
of  $B_j$ containing $xnM$.  If the action of $B_j$
on $M$ is nontrivial then by Proposition~\ref{one-step} there exists 
a unique $m\in M$ such that $xnm$ is contained in a lift of $Q$ in $H$. 
It follows therefore that if the action of $B_j$ on $M$ is 
nontrivial then the set of  $m\in M$ for which $xnm$ is contained
in a lift of $B_j$ is countable. But this is a contradiction
since in fact since every $xnm$, $m\in M$, is contained in a 
one-parameter subgroup of $H$, which necessarily has to be a lift of
some $B_j$ such that $(a,B_j)$ is of type $\cal E$ on $H/M$.  The contradiction shows that $a$ must indeed by
of type~$\cal E$, which proves the theorem. \qed

\begin{remark}
{\rm The criterion in Theorem~\ref{solvable} involves that for 
$a\in A$ there exists a {\it common} one-parameter subgroup $B$ such that for 
every irreducible subquotient $(M,M')$, $M,M'\in \frak N (H)$,  
such that if the action of $a$ on $M/M'$ is trivial the action of $B$ 
is also trivial. It does not suffice to have such a one-parameter for 
each $(M,M')$ individually, and dependent on it. This is illustrated by the following example:  
Let $H$ be the semi-direct product of $\T^2=\{(\rho_1,\rho_2)\mid \rho_1,\rho_2 \in \C, |\rho_1|
=|\rho_2|=1\}$ and $\C^2=\{(z_1,z_2)\mid z_1,z_2\in \C\}$, with the 
conjugation action of $(\rho_1,\rho_2)$ given by $(z_1,z_2)\mapsto
(\rho_1\rho_2z_1, \rho_1\rho_2^{-1}z_2)$. Then $N$ and $A$ as in the earlier
notation may be seen to correspond to $\C^2$ and $\T^2$ respectively, 
(together with the actions involved). The irreducible subquotient 
actions of $A$ on $N$ are seen to correspond to actions on the two copies of 
$\C$, corresponding to the two coordinates. 
Let  $a=(-1,-1) \in \T^2$. Then the action of $a$ on $\C^2$ is trivial.  
We see that there are one-parameter subgroups $B_1$ and $B_2$ 
of $\T^2$ containing $a$ acting trivially on the subquotients corresponding to 
the first
and second coordinates respectively, but there is no one-parameter subgroup 
containing $a$ and acting trivially on both coordinates. Thus the condition 
of Theorem~\ref{solvable} is not satisfied. Hence by the Theorem the 
group is not exponential, as may also be checked directly.

}
\end{remark}

Theorem~\ref{solvable} can be used to deduce the following property, 
known earlier (\cite{W}, Corollary~3.18). 

\begin{corollary}\label{cor}
Let $H$ be an exponential solvable Lie group. Then the center of 
$H$ is connected. 
\end{corollary} 

\proof Let $N$ be the nilradical of $H$ and $A=H/N$. Let
$z$ be an element contained in the center of $H$ and $a=zN\in A$. Since 
$H$ is exponential $xn$ is exponential for all $n\in N$. Also, $z$ being contained
in the center,  the action of $a$ 
on any irreducible subquotient is trivial. Hence by Theorem~\ref{solvable} 
there exists a one-parameter subgroup $B$ of $A$ containing $a$ such that 
the action of $B$ on any irreducible subquotient is trivial. Let $L$ be the subgroup of 
$H$ containing $N$ and such that $B=L/N$. Then the preceding 
observation   implies that 
$L$ is a connected nilpotent Lie group. Since $N$ is the nilradical we get that $L=N$. 
Therefore $z$ is contained in the center of $ N$, say $Z$. Now $Z$ is a connected abelian 
Lie group, and the center of $H$ is precisely the set of fixed points of the action of $A$ on $Z$, 
induced by the $A$-action on $N$. It can be seen that as $A$ is connected, the set of
fixed points is a connected subgroup of $Z$. Hence the  center of $H$ is connected. \qed

We also deduce the following characterization of exponential Lie group, 
which is in a sense the main nontechnical part in the characterization of exponentiality
of solvable Lie groups in  
Theorem~3.17 of \cite{W}. We follow the terminology as in \cite{W};
however the symbols chosen are different, to avoid clash with the
notation in the rest of this paper. 

\begin{corollary}
Let $H$ be a connected solvable Lie group, $\frak H$ be the Lie algebra 
of $H$,  and $C$ be a Cartan subgroup of $H$. 
Then $H$ is an exponential Lie group if and only if for all nilpotent
 elements $\nu \in \frak H$ the centralizer 
of $\nu$ in $C$, namely $Z_C(\nu):=\{y\in C\mid \Ad (y) (\nu)=\nu\}$,   is 
connected. 
\end{corollary} 

\proof Let $N$ be the nilradical of $H$ and $\frak N$ be the Lie subalgebra
corresponding to $N$. Then we have $H=CN$.  As $\frak H$ is solvable, all 
nilpotent elements in $\frak H$
are contained in  $\frak N$. Also, $\frak N$ has a decomposition
as $\oplus_{s\in S}\frak N_s$, with $S$ an indexing set, such that 
the following holds: for each $s\in S$, 
$\frak N_s$ is
$\Ad C$-invariant, for each $y\in C$ the  eigenvalues 
of the restriction of $\Ad y$ to  
$\frak N_s$ consist of a complex conjugate pair, say $\lambda (y,s)$ and 
$\overline {\lambda (y,s)}$ (only one when real), and there exists an order 
on $S$,
denoted by $\geq$, such that
for any $t\in S$, $\sum_{s\geq t}\frak N_s$ 
is a Lie ideal in $\frak H$; more precise relations can be written down for 
commutators of the $\frak N_s$'s but we do not need that here. We note
also that $C$ is a connected nilpotent Lie group and hence given $y\in C$
and a one-parameter subgroup $B'$ of $CN/N$ containing $yN$ there exists a 
one-parameter $B$ of $C$ containing $y$ such that $B'=BN/N$. 

Now suppose that $H$ is exponential and let a nilpotent $\nu$ element be
given. Then we have $\nu =\sum_{s\in S}\nu_s$, with $\nu_s\in \frak N_s$ for
all $s\in S$. Let  $y\in C$
and $a=yN/N$. By Theorem~\ref{solvable}, together with one of the 
observations 
above, there exists a one-parameter 
subgroup $B$ of $C$ containing $y$  such that for any irreducible 
subquotient $(M,M')$ such that  
the action of $y$ on $M/M'$ is trivial, the action of $B$ on $M/M'$ is also 
trivial. 
Now, $\Ad y (\nu)=\nu$ if and only if $\lambda (y,s)=1$ for all $s$ 
such that $\nu_s\neq 0$. Consider any  $t\in S$ such that $\lambda (y,t)=1$. 
Let $M$ be the closed subgroup with Lie algebra   $\sum_{s\geq t}\frak N_s$ 
and $M'$ be a closed connected normal subgroup of $H$, properly contained 
in $M$, such that the Lie algebra of $M'$ contains $\frak N_s$, for all $s> t$ (that 
is $s\geq t$ and $s\neq t$), and $(M,M')$  is an irreducible quotient.  
As $\lambda (y,t)=1$ it follows that 
the action of $y$ on $M/M'$ is trivial. Hence the action of $B$ on $M/M'$
is also trivial and in turn  $\lambda (b , t)=1$ for all $b\in B$. 
Thus we get that $\Ad b (\nu)=\nu$ for all $b\in B$. This shows that 
the centralizer of $\nu$ in $C$ is connected. 
 
Now suppose that $Z_C(\nu)$ is connected for all 
$\nu \in \frak N$. Let $x\in H$ be given and let $a=xN\in H/N$. Since 
$H=CN$ there there exists  $y\in C$ such that $a=yN/N$.  Let $S'=\{s\in S\mid
\lambda (y,s)=1\}$. There exists an element $\nu=\sum_{s\in S'}\nu_s$, such 
that each $\nu_s$, $s\in S'$, is a nonzero element of $\frak N_s$ such 
that $\Ad y(\nu_s)=\nu_s$. In particular $\nu$ is fixed by $\Ad y$, and 
since   $Z_C(\nu)$ is connected and nilpotent we get that there exists a one-parameer 
subgroup $B$ of $C$ containing $y$ such that $\Ad b (\nu)=\nu$ for all $b\in B$. 
Hence $\Ad b (\nu_s)=\nu_s$,  and in turn  
$\lambda (b,s)=1$, for all $s\in S'$ and $b\in B$.
Now let $B'=BN/N$. Then $B'$  is a one-parameter subgroup of $H/N$ 
containing $a$; we shall
show that $(a,B')$ is of type $\cal E$.  

Let $(M,M')$ be 
an irreducible subquotient such that the action of $a$ on $M/M'$ is trivial. 
Then there exists $t\in S'$ such that $M$ is contained in 
$\sum_{s\geq t}\frak N_s$ and $M'$ contains $\sum_{s> t}\frak N_s$. Since 
$\lambda (b,s)=1$ for all $s\in S'$ and $b\in B$, it follows that
the action of $B'=BN/N$ on  $M/M'$ is trivial. Thus  $(a,B')$ is of type 
$\cal E$ and hence by Theorem~\ref{solvable}
$x$ is exponential in $H$. Therefore $H$ is an exponential Lie group. \qed

\section{Proof of Theorem~\ref{thm}}

We shall now deduce Theorem~\ref{thm}; we follow the notations  $G$, $R$, and $S$  as in the statement
 of the theorem in \S\,1. Since $S$ satisfies condition $\frak U$ there 
exists a unipotent one-parameter subgroup $U$ of $S$ such that the centralizer 
of $U$ in $S$ contains no compact subgroup of positive dimension. We note the
following. 

\begin{lemma}\label{lem}
Let $S$ and $U$  be as above. Let $u$ be a nontrivial element of $U$. Then $U$ is 
the only one-parameter subgroup of $S$ containing $u$. 
\end{lemma}  

\proof We may assume $U=\{u_t\}$ and $u=u_1$. It suffices to show that if  $\{x_t\}$ is any one-parameter 
subgroup of $S$ such that $x_1=u$ then $x_t=u_t$ for all $t$. Let $s\in \R$ be arbitrary. 
Since $x_s$ commutes
with $x_1=u$ and $\{u_t\}$ is a unipotent one-parameter subgroup it follows that  $x_s$
commutes with $u_t$ for all $t\in \R$. Since this holds for all $s\in \R$ it follows that 
$\{x_{-t}u_t\} $ is a one-parameter subgroup. Moreover, since $x_{-1}u_1=u^{-1}u =e$,
the identity element, $\{x_{-t}u_t\}$ is a compact subgroup. Since it is contained in 
the centraliser of $U$, the condition on $U$ implies that the subgroup is trivial, namely
$x_t=u_t$ for all $t$. This proves the Lemma. \qed

Now let  
$\eta: G\to S$ 
be the canonical quotient map of $G$ onto $S$, and  let $H=\eta^{-1}(U)$. 
We note that $U$ is a closed connected subgroup of $S$ and hence $H$ 
is a connected Lie group.  
Since $N$ is the nilradical of $G$, $G/N$ is reductive, 
and in particular it follows that $H/N$ is abelian; in particular $H$ is a solvable Lie
group. 

We now first prove the following. 

\begin{proposition}\label{H-exp}
Any  $h\in H$ which is not contained in $R$ is exponential in $H$. 
\end{proposition}

\proof  Let $h\in H\backslash R$ be given. Since $G$ is 
exponential there exists a one-parameter subgroup $P$ of $G$  
containing  $h$. 
Consider the one-parameter subgroup $\eta (P)$. 
We note that   $\eta (h)$ is a nontrivial element of $U$. Hence  by Lemma~\ref{lem} we 
have  $\eta (P)=U$. Thus $P$ is contained in $\eta^{-1}(U)=H$. Therefore  
$h$ is exponential in~$H$.  \qed

\medskip
\noindent{\it Proof of Theorem~\ref{thm}}: 
We first note that in proving the theorem $N$ may be assumed to be simply connected: 
Let $C$ be maximal compact subgroup of $N$. Then $C$ is a connected subgroup 
contained in the centre
of $G$. We note that  the hypothesis of the theorem holds 
for $G/C$, and upholding that the desired conclusion for $G/C$ (namely showing 
that $R/C$ is exponential) implies that $R$ is exponential, as desired. We 
may therefore assume without loss of generality that $C$ is trivial. 
Equivalently this means that $N$ is simply connected. 

Now let $x\in R$ be given. Let $H$ be
the solvable group introduced above and $y=ux \in H$, where $u$ is a nontrivial element of $U$. 
We note that since $H/R$
is topologically isomorphic to $\R$, to prove that $x$ is exponential in $R$ it suffices
to prove that it is exponential in $H$. Now
let $A=H/N$ and $a, a' $ 
be the elements $a=xN$ and $a'=yN$. 
By Proposition~\ref{H-exp} 
$yn$ is exponential in $H$ for all $n\in N$. Hence by Theorem~\ref{solvable}, 
there exists a one-parameter subgroup $B'$ of $A$ containing $a'$ such that 
 for any irreducible subquotient 
$(M,M')$, $M, M'\in \frak N(H)$, for which the action of $a'$ on $M/M'$ is trivial, 
the action of $B'$  on $M/M'$ is also trivial. Let $B'=\{b'_t\}$ and $U=\{u_t\}$, with $b'_1=a'$ and 
$u_1=u$. Let $B=\{b_t\}$ be the one-parameter subgroup of $A$ defined by $b_t=(u_{-t}N)b'_t$
for all $t\in R$; in particular $b_1=u^{-1}a'=a$.  
We note that as $U$ is a unipotent one-parameter subgroup of $S$, the action of $U$ 
on any irreducible subquotient $M/M'$ as above  is trivial. Hence on any 
irreducible subquotient the actions of $b_t$ and $b'_t$ coincide for each $t\in \R$. 
Therefore  we get that for any irreducible subquotient 
$(M,M')$, $M, M'\in \frak N(H)$, for which the action of $a$ on $M/M'$ is trivial, 
the action of $B$  on $M/M'$ is trivial. Thus condition~(i)
of Theorem~\ref{solvable} holds for $x$, and the theorem implies that $xn$ is 
exponential for all $n\in N$; in particular $x$ is exponential in $H$, and hence, as seen above, 
also in $R$. This proves the theorem. \qed

\begin{corollary}
Let $G$ be a Lie group as in Theorem~\ref{thm} and  $R$ be its radical.  Let $Z(G)$ 
and $Z(R)$ be the centers of $G$ and $R$ respectively. Then $Z(R)$  and $Z(G)\cap R$ 
are connected. 

\end{corollary}

\proof By Corollary~\ref{cor} $Z (R)$ is connected. Clearly $Z(G)\cap R$ is the set of fixed 
points of the canonical action of $G/R$ on $Z(R)$, induced by the conjugation action of $G$. As $Z(R)$ is abelian and $G/R$ is connected it follows that the set of fixed points is a connected subgroup. \qed 

\section{Groups satisfying condition $\frak U$}

In this section we consider semisimple Lie groups $S$ satisfying condition $\frak U$ and 
prove Theorem~\ref{ss}. We 
shall follow the notation as in \S\,1; we recall in particular that $M$ denotes the subgroup 
consisting of all elements in a maximal compact subgroup $K$ of $S$, commuting with all 
elements in a  maximal subgroup $\Ad$-diagonalisable over $\R$.

\begin{remark}\label{semireg}
{\rm  We note that when $M$ is abelian there exists a unipotent one-parameter subgroup 
whose centraliser is contained in the center of $S$. It suffices
to see this in the case when $S$ has trivial center, so $S$ may be taken to be the group
of $\R$-elements of a semisimple algebraic group defined over $\R$. In this case $M$ being abelian 
is equivalent to the group being quasi-split, namely that there exists a Borel subgroup 
defined over $\R$.  Under that condition it follows from  \cite{R}, Proposition~5.1 (see also \cite{C} for  a generalisation of the result to non-semisimple 
algebraic groups) and \cite{A} and  Corollary~5 (page 114) that there exists $u\in N$ such that
the centraliser of $u$ in $S$ is contained in $P$. On the other  
the conjugation action of $M$ action on $N$ has no fixed points, so the centraliser of $u$, and hence
of any one-parameter subgroup $U$ containing $u$, has no compact subgroup of positive dimension. 
Hence $S$ satisfies condition $\frak U$. 
}
\end{remark}

\medskip
\noindent {\it Proof of Theorem~\ref{ss}}:  Since $M$ is abelian by Remark~\ref{semireg} there exists a unipotent one-parameter subgroup, say  $U$, in $S$ whose
centraliser is contained in the (discrete) center of $S$; in particular the centraliser of 
$U$ in $S$ contains no subgroup of positive dimension. Hence condition~$\frak U$ is satisfied for $S$.

For split semisimple Lie groups the subgroup $M$ as above is trivial, and for complex Lie groups 
it is abelian. It therefore follows from the above results the condition $\frak U$ holds in 
these cases. \qed

\vskip15mm

\begin{flushleft}
S.G. Dani\\
Department of Mathematics\\
Indian Institute of Technology Bombay\\
Powai, Mumbai 400076\\
India

\medskip
E-mail: {\tt sdani@math.iitb.ac.in}
\end{flushleft}

\end{document}